\documentclass[12pt,reqno,a4paper]{amsart}
\usepackage{mathrsfs}
\usepackage{pifont}
\usepackage{latexsym,amsmath,amsfonts,amssymb,amsthm,epic,graphicx}
\theoremstyle{plain}
\newtheorem{Thm}{Theorem}
\newtheorem{Lem}{Lemma}

\theoremstyle{definition}
\newtheorem*{Ack}{Acknowledgment}
\theoremstyle{remark}

\def\Z{\mathbb Z}

\def\N{\mathbb N}

\def\A{\mathcal A}
\def\B{\mathcal B}
\def\K{\overrightarrow{K}}
\def\G{\overrightarrow{G}}
\def\HH{\mathscr{H}}
\def\1{\mathbf 1}

\def\pmod #1{\ ({\rm{mod}}\ #1)}

\begin{document}
\title{Clique numbers of graphs and irreducible exact $m$-covers of $\Z$}
\date{}
\author{Hao Pan}
\address{Department of Mathematics, Shanghai Jiaotong University, Shanghai
200240, People's Republic of China} \email{haopan79@yahoo.com.cn}
\author{Li-Lu Zhao}
\address{Department of Mathematics, Nanjing University\\
 Nanjing 210093, People's Republic of China} \email{zhaolilu@gmail.com}
\begin{abstract}
For each $m\geq 1$, we construct a graph $G=(V,E)$ with
$\omega(G)=m$ such that
$$
\max_{1\leq i\leq k}\omega(G[V_i])=m
$$
for arbitrary partition $\{V_1,\ldots,V_k\}$ of $V$, where
$\omega(G)$ is the clique number of $G$ and $G[V_i]$ is the induced
graph of $G$ with the vertex set $V_i$. Using this result, we show
that for each $m\geq 2$ there exists an exact $m$-cover of $\Z$
which is not the union of two 1-covers.
\end{abstract} 
\thanks{}
\maketitle

\section{Introduction}
\setcounter{equation}{0} \setcounter{Thm}{0} \setcounter{Lem}{0}
\setcounter{Cor}{0}

In his proof of the existence of irreducible exact $m$-covers of
$\Z$ (the notions will be introduced soon), Zhang proved the
following graph-theoretic result \cite[Lemma 2]{Zhang1991}:
\begin{Thm}
\label{t1} For every $m\geq 1$, there exists a graph $G=(V,E)$
satisfying the following properties:

\medskip
$\omega(G)=m$, where $\omega(G)$ is the clique number of $G$, i.e.,
the maximal order of the complete subgraphs of $G$ . And if the
vertex set $V$ is arbitrarily split into two non-empty subsets $V_1$
and $V_2$, then
$$
\omega(G[V_1])+\omega(G[V_2])>\omega(G),
$$
where $G[V_i]$ denotes the induced subgraph of $G$ with the vertex
set $V_i$.
\end{Thm}

In this paper, our main purpose is to give an extension of Zhang's
result as follows:
\begin{Thm}
\label{t2} For every $m\geq 1$ and $k\geq 2$, there exists a graph
$G=(V,E)$ with $\omega(G)=m$ satisfying the following property:

\medskip
If the vertex set $V$ is arbitrarily split into $k$ subsets $V_1,
V_2,\ldots,V_k$, then
$$
\max_{1\leq i\leq k}\omega(G[V_i])=\omega(G).
$$
\end{Thm}

For an integer $a$ and a positive integer $n$, let $a(n)$ denote the
residue class $\{x\in\Z:\, x\equiv a\pmod{n}\}$. For a finite system
$\A=\{a_t(n_t)\}_{t=1}^s$, define the covering function $w_\A$ over
$\Z$ by
$$
w_\A(x):=|\{1\leq t\leq s:\, x\in a_t(n_t)\}|.
$$
If $w_\A(x)\geq m$ for each $x\in\Z$, we say that a system $\A$ is
an {\it $m$-cover} of $\Z$. In particular, we call $\A$ an {\it
exact $m$-cover} provided that $w_\A(x)=m$ for all $x\in\Z$. The
covers of $\Z$ was firstly introduced by Erd\H os \cite{Erdos50} and
has been investigated in many papers (e.g.,
\cite{Newman71,Porubsky75,Znam75,Simpson85,BergerFelzenbaumFraenkel88,Sun95,Sun96,Sun00,
Chen03,GuoSun05}).

Suppose that $\A_1$ is an $m_1$-cover and $\A_2$ is an $m_2$-cover,
then clearly $\A=\A_1\cup\A_2$ forms an $(m_1+m_2)$-cover.
Conversely, Porubsk\'y \cite{Porubsky76} asked whether for each
$m\geq 2$ there exists an exact $m$-cover of $\Z$ which cannot be
split into an exact $n$-cover and an exact $(n-m)$-cover with $1\leq
n<m$. Choi gave such a example for $m=2$:
$$
\A=\{1(2);0(3);2(6);0,4,6,8(10);1,2,4,7,10,13(15);5,11,12,22,23,29(30)\}.
$$
In \cite{Zhang1991}, using Theorem \ref{t1}, Zhang gave an
affirmative answer to Porubsk\'y's problem. This shows that the
results on $m$-covers of $\Z$ is essential. In \cite{Sun}, Sun
established a connection between $m$-covers of $\Z$ and zero-sum
problems in abelian $p$-groups. For more related results, the
readers may refer to \cite{Sun92,Sun97,Sun99}

On the other hand, for each $m\geq 2$, Pan and Sun \cite[Example
1.1]{PanSun07} constructed an $m$-cover of $\Z$ (though not exact)
which even is not the union of two $1$-covers! As an application of
Theorem \ref{t2}, we have a common extension of the above two
results:
\begin{Thm}
\label{cover} For each $m\geq 2$, there exists an exact $m$-cover of
$\Z$ which is not the union of two $1$-covers.
\end{Thm}

We shall prove Theorem \ref{t2} in the next section, and the proof
of Theorem \ref{cover} will be given in Section \ref{s3}.

\section{Proof of Theorem \ref{t2}}
\label{s2} \setcounter{equation}{0} \setcounter{Thm}{0}
\setcounter{Lem}{0} \setcounter{Cor}{0}

\begin{Lem}
\label{Di} Suppose that $G=(V,E)$ is a connected simple graph and
$v_0$ is a vertex of $G$. Then there exists an oriented graph $\G$
arising from $G$, which satisfies that:

\noindent(i)\quad $\G$ doesn't contains any directed cycle.

\noindent(ii)\quad For any vertex $u\in V\setminus\{v_0\}$, there
exists a directed path of $\G$ from $v_0$ to $u$.
\end{Lem}
\begin{proof}
We use induction on $|V|$. There is nothing to do when $|V|=1$ or
$2$. Now assume that $|V|>0$ and our assertion holds for any smaller
value of $|V|$. Let $V'=V\setminus\{v_0\}$ and $G'=G[V']$. Suppose
that $v_1,\ldots,v_s\in V'$ are all vertex adjacent to $v_0$ in $G$.
By the induction hypothesis, there exists an oriented graph
$\overrightarrow{G'}$ obtained from $G'$, satisfying the properties
(i) and (ii) for the vertex $v_1$. Now we direct the edge $v_0v_i$
from $v_0$ to $v_i$ for $1\leq i\leq k$, and preserve the direction
of each edge in $\overrightarrow{G'}$. Thus we obtain an oriented
graph $\G$. Clearly $\G$ doesn't contain any directed cycle since
$v_0$ can't lie in any directed cycle. And for any $u\in
V\setminus\{v_0,v_1\}$, since there exists a directed path of
$\overrightarrow{G'}$ from $v_1$ to $u$, the property (ii) is also
satisfied.
\end{proof}

\begin{Lem}
\label{My} For every $k\geq 1$, we can construct a $k$-chromatic
graph without any triangle.
\end{Lem}
\begin{proof} The reader may refer to \cite{Mycielsky55} (or \cite[Chapter
5, Exercise 23]{Diestel05}) for the construction of such graph. In
fact, with help of his probabilistic method, Erd\H os \cite{Erdos59}
proved that there exist the graphs having arbitrarily large girths
and chromatic numbers. \end{proof}

\begin{proof}[Proof of Theorem \ref{t2}]

Let $K=(V_K,E_K)$ be a $(k+1)$-chromatic graph without any triangle.
Let $u_0$ be a vertex of $K$. Then there exists an oriented graph
$\K$ arising from $K$, which satisfies the properties (i) and (ii)
of Lemma \ref{Di} for the vertex $u_0$. Let $n=|V_K|$ and suppose
that $u_0,u_1,\ldots,u_{n-1}$ are all vertices of $K$. For $1\leq
i\leq n-1$, let $l_i$ denote the length of the longest directed path
from $u_0$ to $u_i$ in $\K$. By the property (ii) of Lemma \ref{Di},
these $l_i$ are well-defined. Let $l=\max_{1\leq i\leq n-1} l_i$,
and for $1\leq j\leq l$ let
$$
D_j=\{1\leq i\leq n-1;\,l_i=j\}
$$
In particular, we set $D_0=\{0\}$. For $1\leq i\leq n-1$, let
$$
A_i=\{0\leq i'\leq n-1:\,\overrightarrow{u_{i'}u_i}\text{ lies in
}\K\},
$$
where we denote by $\overrightarrow{xy}$ the directed edge from $x$
to $y$. In particular, we set $A_0=\emptyset$.
\begin{Lem}
\label{Ai} For $1\leq j\leq l$, we have
\begin{equation}
\bigcup_{u_i\in D_j}A_i\subseteq\bigcup_{0\leq j'\leq j-1} D_{j'}.
\end{equation}
\end{Lem}
\begin{proof}
Assume on the contrary that there exist $u_i\in D_j$ and $i'\in A_i$
such that $u_{i'}\not\in\bigcup_{0\leq j'\leq j-1} D_{j'}$. From the
definition of $D_{j'}$, we know that there exists a path from $u_0$
to $u_{i'}$ with the length at least $j$. If $u_i$ doesn't lie in
this path, then we get a path from $u_0$ to $u_{i}$ with the length
at least $j+1$, since the direction of the edge
$\overrightarrow{u_{i'}u_i}$ is from $u_{i'}$ to $u_i$. On the other
hand, if $u_i$ lies in this path, then clearly we get a directed
cycle from $u_{i}$ to $u_{i'}$, next to $u_i$. This also leads to a
contradiction with the property (i) of Lemma \ref{Di}.
\end{proof}

\begin{Lem}
\label{Dj} $D_j\not=\emptyset$ for each $1\leq j\leq l$.
\end{Lem}
\begin{proof}
Clearly $D_l\not=\emptyset$. Let $u_{i_l}$ be a vertex in $D_l$.
Then there exists a directed path in $\K$ from $u_0$ to $u_{i_l}$
with the length $l$. Suppose that this path is
$$
u_0\to u_{i_1}\to u_{i_2}\to\cdots\to u_{i_{l-1}}\to u_{i_l}.
$$
We claim that $i_j\in D_j$ for each $1\leq j\leq l$. We use
induction on $j$. Clearly our assertion holds when $j=l$. Assume
that $j<l$ and $i_{j+1}\in D_{j+1}$. Clearly $l_{i_j}\geq j$ since
$u_0\to u_{i_1}\to\cdots\to u_{i_j}$ is a directed path with the
length $j$. On the other hand, by Lemma \ref{Ai}, we have
$$
u_{i_j}\in A_{i_{j+1}}\subseteq\bigcup_{0\leq j'\leq j} D_{j'}.
$$
Hence $l_{i_j}\leq j$. So $l_{i_j}=j$ and $i_j\in D_j$. We are done.
\end{proof}

We shall use induction on $m$ to prove Theorem \ref{t2}. The case
$m=1$ is trivial. Now assume that $m\geq2$ and our assertion holds
for $m-1$. That is, there exists a graph
$G^{(m-1)}=(V^{(m-1)},E^{(m-1)})$ with $\omega(G^{(m-1)})=m-1$
satisfying that
$$
\max_{1\leq i\leq k}\omega(G^{(m-1)}[V_i])=m-1
$$
for arbitrary partition $V_1,\ldots,V_k$ of $V^{(m-1)}$.

First, we shall create $n$ graphs $H_0,H_1,\ldots,H_{n-1}$. $H_0$ is
a graph only having a vertex $x_0$. For each $i\in D_{1}$, $H_i$ is
one copy of $G^{(m-1)}$. Similarly, for $2\leq j\leq l$ and every
$i\in D_{j}$, assuming $H_{i'}$ have been created for all
$i'\in\bigcup_{0\leq j\leq j-1} D_{j'}$, let $H_{i}$ be
$$
h_i:=\prod_{i'\in A_i}|V(H_{i'})|
$$ disjoint copies
of $G^{(m-1)}$, where $V(H_{i'})$ denotes the vertex set of
$H_{i'}$.

Next, we shall add some edges between the vertices of $H_i$ and the
vertices of $H_{i'}$, for $0\leq j<j'\leq l$, $i\in D_{j}$ and
$i'\in D_{j'}$. For every $i\in D_{1}$, we join $x_0$ and $H_i$,
i.e., join $x_0$ and all vertices of $H_i$. Below we shall
inductively add the edges incident with the vertices of $H_i$ for
every $2\leq j\leq l$ and $i\in D_j$. Suppose that $2\leq j\leq d$,
$i\in D_j$ and $A_i=\{i_1',\ldots,i_s'\}$ with $i_1'<\cdots<i_s'$.
Assume that we have added the edges between the vertices of
$H_{i_1'}$ and $H_{i_2'}$, for every $0\leq j_1'<j_2'\leq j-1$ and
$i_1'\in D_{j_1'},\ i_2'\in D_{j_2'}$. Recall that $H_i$ is formed
by $h_i$ disjoint copies of $G^{(m-1)}$. Let $\psi_i$ be an
arbitrary $1-1$ projection from $V(H_{i_1'})\times\cdots\times
V(H_{i_s'})$ to $\{1,\ldots,h_i\}$, where
$V(H_{i_1'})\times\cdots\times V(H_{i_s'})$ denotes the Cartesian
product of $H_{i_1'},\ldots,H_{i_s'}$. Then for each
$(w_1,\ldots,w_s)\in V(H_{i_1'})\times\cdots\times V(H_{i_s'})$, we
join the vertices $w_1,\ldots,w_s$ to the
$\psi_i(w_1,\ldots,w_s)$-th copy of $G^{(m-1)}$ in $H_i$. Taking the
above processes from $j=2$ to $l$, we obtain the desired graph
$G^{(m)}=(V^{(m)},E^{(m)})$.

The remainder task is to show that $G_m$ certainly satisfies our
requirements. Clearly $\omega(G^{(m)})\geq m$ since
$\omega(G^{(m-1)})=m-1$ and $x_0$ is adjacent to all vertices of at
least one copy of $G^{(m-1)}$. Let $\Omega$ be an arbitrary complete
subgraph of $G^{(m)}$. We need to prove that $\Omega$ has at most
$m$ vertices. Let $U_i$ be the set of all vertices of $\Omega$ lying
in $H_i$. Notice that for distinct $i$ and $i'$, if there exist
$w\in H_{i}$ and $w'\in H_{i'}$ such that $ww'\in E^{(m)}$, then
either $i\in A_{i'}$ or $i'\in A_{i}$, i.e., $u_i$ and $u_{i'}$ are
adjacent in the graph $K$. Since $K$ doesn't contain any triangle,
we have $|\{i:\, U_i\not=\emptyset\}|\leq 2$. There is noting to do
if $\Omega$ is completely contained in one $H_i$, since
$\omega(H_i)=\omega(G^{(m-1)})=m-1$. Suppose that there exist
distinct $i,i'$ such that $U_i, U_{i'}\not=\emptyset$. Without loss
of generality, assume that $i'\in A_i$. Observe that distinct
vertices of $H_{i'}$ are joint to distinct copies of $G^{(m-1)}$ in
$H_i$. So we must have $|U_{i'}|=1$. Hence
$$
|V(\Omega)|=|U_{i}|+|U_{i'}|\leq\omega(G^{(m-1)})+1=m.
$$

Now assume that the vertex set $V^{(m)}$ is split into $k$ disjoint
subsets $V_1,\ldots,V_k$. Without loss of generality, we may assume
that $x_0\in V_1$. Let $U_{i,g}^{(t)}$ be the set of the common
vertices of $V_t$ and the $g$-th copies of $G^{(m-1)}$ in $H_i$. By
the induction hypothesis, we know that
$$
\max_{1\leq t\leq
k}\omega(G^{(m)}[U_{i,g}^{(t)}])=\omega(G^{(m-1)})=m-1
$$
for every $1\leq i\leq n$ and $1\leq t\leq h_i$. For every $i\in
D_1$, let $g_i=1$,
$$
t_{i}=\min\{1\leq t\leq k:\,\omega(G^{(m)}[U_{i,1}^{(t)}])=m-1\}
$$
and arbitrarily choose a vertex $w_i\in U_{i,1}^{(t_i)}$. Below we
shall determine $g_i$, $t_i$, $w_i$ inductively for $2\leq j\leq l$
and $i\in D_j$. Assume that $j\geq 2$ and we have determined $g_i$,
$t_i$, $w_i$ for all
$$
i\in\bigcup_{1\leq j'\leq j-1} D_{j'}.
$$
Then for $i\in D_j$, supposing $A_i=\{i_1',\ldots,i_s'\}$ with
$i_1'<\cdots<i_s'$, let $g_i=\psi_i(w_{i_1'},\ldots,w_{i_s'})$,
$$
t_{i}=\min\{1\leq t\leq k:\,\omega(G^{(m)}[U_{i,g_i}^{(t)}])=m-1\}
$$
and let $w_i$ be an arbitrary vertex in $U_{i,g_i}^{(t_i)}$. In
particular, we set $t_0=1$ and $w_0=x_0$.

Now we shall color the vertices of $K$ with $k$ colors. For $0\leq
i\leq n-1$, let the vertex $u_i$ be colored with the $t_i$-th color.
Since $K$ is not $k$-colorable, there exist distinct $0\leq i,i'\leq
n-1$ such that $t_i=t_{i'}$ and $u_iu_{i'}\in E_K$, i.e., either
$i\in A_{i'}$ or $i'\in A_{i}$. Without loss of generality, assume
that $i'\in A_i$. Notice that $w_{i'}\in U_{i',g_{i'}}^{(t_i)}$ and
$w_{i'}$ is adjacent to all vertices of the $g_i$-th copies of
$H_i$. Also, we have $G^{(m)}[U_{i,g_i}^{(t_i)}]$ contains an
$(m-1)$-complete subgraph. Thus we get an $m$-complete subgraph of
$G^{(m)}[U_{i,g_i}^{(t_i)}\cup\{w_{i'}\}]$, which is also a subgraph
of $G^{(m)}[V_{t_i}]$. We are done.

\end{proof}

\section{Proof of Theorem \ref{cover}}
\label{s3} \setcounter{equation}{0} \setcounter{Thm}{0}
\setcounter{Lem}{0} \setcounter{Cor}{0}

For a system $\A=\{a_t(n_t)\}_{t=1}^s$ and a graph $G=(V,E)$ with
$V=\{v_1,\ldots,v_s\}$, we say $G$ is an intersection graph of $\A$
if
$$
a_i(n_i)\cap a_j(n_j)\not=\emptyset\Longleftrightarrow\text{the edge
 }v_iv_j\in E
$$
for any $1\leq i<j\leq s$. The following result \cite[Theorem
1]{Zhang1991} is due to Zhang, although we give a slightly different
proof here for the sake of completeness.
\begin{Lem}
\label{l31} For each graph $G=(V,E)$ with $|V|=s$, there exists a
system $\A=\{a_t(n_t)\}_{t=1}^s$ such that $G$ is an intersection
graph of $\A$.
\end{Lem}
\begin{proof} We use induction on $s$. The cases $s=1$ and $s=2$ are
trivial. Assume that $s>2$ and our assertion holds for $s-1$.
Suppose that $V=\{v_1,\ldots,v_s\}$. Let $V'=V\setminus\{v_s\}$ and
$G'=G[V']$. Let $\A'=\{a_t'(n_t')\}_{t=1}^{s-1}$ be a system such
that $G'$ is an intersection graph of $\A'$. Let
$p_1,\ldots,p_{s-1}$ be some distinct primes greater than
$\max\{n_1',\ldots,n_{s-1}'\}$. For each $1\leq t\leq s-1$, let
$n_t=n_t'p_t$ and $a_t$ be an integer such that $a_t\equiv
a_t'\pmod{n_t'}$ and $a_t\equiv 1\pmod{p_t}$. Let $n_s=p_1\cdots
p_{s-1}$ and $a_s$ be an integer such that
$$
a_s\equiv\begin{cases} 1\pmod{p_t}&\text{if the
edge }v_tv_s\in E,\\
0\pmod{p_t}&\text{if the edge }v_tv_s\not\in E
\end{cases}
$$
for $1\leq t\leq s-1$. Since $a_i(n_i)\cap a_j(n_j)\not=\emptyset$
if and only if $(n_i,n_j)\mid a_i-a_j$, it is easy to see that $G$
is an intersection graph of the system $\A=\{a_t(n_t)\}_{t=1}^s$.
\end{proof}

Suppose that $G=(V,E)$ is an intersection graph of
$\A=\{a_t(n_t)\}_{t=1}^s$. By the Chinese remainder theorem, for a
subset $I\subseteq\{1,\ldots,k\}$, if $a_i(n_i)\cap
a_j(n_j)\not=\emptyset$ for any $i,j\in I$, then $\bigcap_{i\in
I}a_i(n_i)\not=\emptyset$. Hence we have
$$
\omega(G)=\max\{w_\A(x):\,x\in\Z\},
$$
by recalling that $w_\A(x)=|\{1\leq i\leq s:\, x\in a_s(n_s)\}|$.

\begin{proof}[Proof of Theorem \ref{cover}]
Let $G=(V,E)$ be the graph satisfying the properties in Theorem
\ref{t2} for $k=2$. Assume that $|V|=s$. By Lemma \ref{l31}, there
exists a system $\A=\{a_t(n_t)\}_{t=1}^s$ such that $G$ is an
intersection graph of $\A$. We claim that for any partition
$\{\A_1,\A_2\}$ of $\A$,
$$
\max_{i=1,2}\omega_{\A_i}=\omega_\A,
$$
where
$$
\omega_\A=\max\{w_\A(x):\, x\in\Z\}.
$$
In fact, letting $V_i\subseteq V$ be the set of vertices concerning
those arithmetic progressions in $\A_i$, we have $G[V_i]$ is an
intersection graph of $\A_i$. Hence
$$
\max_{i=1,2}\omega_{\A_i}=
\max_{i=1,2}\omega(G[V_i])=\omega(G)=\omega_\A.
$$
Since $\omega(G)=m$, $w_\A(x)\leq m$ for every $x\in\Z$. So we may
choose integers $b_1,\ldots,b_r$ such that
$\B=\A\cup\{b_j(N)\}_{j=1}^r$ forms an exact $m$-cover, where $N$ is
the least common multiple of $n_1,\ldots,n_s$. If $\B$ is
arbitrarily split into $\B_1$ and $\B_2$, then
$$
\max_{i=1,2}\omega_{\B_i}\geq\max_{i=1,2}\omega_{\B_i\cap\A}=\omega_{\A}=\omega_{\B}.
$$
Hence there exists an integer $x$ such that $w_{\B_1}(x)=m$ or
$w_{\B_2}(x)=m$. Without loss of generality, assume that
$w_{\B_1}(x)=m$. Then $w_{\B_2}(x)=w_{\B}(x)-w_{\B_1}(x)=0$, whence
$\B_2$ is not a $1$-cover.
\end{proof}

\section{A Further Remark}
\label{s4} \setcounter{equation}{0} \setcounter{Thm}{0}
\setcounter{Lem}{0} \setcounter{Cor}{0}

We may consider a general problem. Let $\HH$ be a set of graphs such
that for any $G\in\HH$, all induced subgraphs of $G$ are also
contained in $\HH$. Suppose that $\psi$ be a projection from $\HH$
to $\N=\{0,1,2,\ldots\}$. We may ask whether for every $m\geq 0$ and
$k\geq 2$, there exists a graph $G=(V,E)\in\HH$ with $\psi(G)=m$
satisfying that
$$
\psi(G)\in\{\psi(G[V_1]),\psi(G[V_2]),\psi(G[V_k])\}
$$
for any $k$-partition $\{V_1,V_2,\ldots,V_k\}$ of the vertex set
$V$.

Let $l(G)$ denote the length of the longest path of $G$. Then we
have the following negative result for $l(\cdot)$.
\begin{Thm}
Let $G=(V,E)$ be a graph having at least one edge. Then there exists
a partition $\{V_1,V_2\}$ of the vertex set $V$ such that
$$
l(G[V_1])<l(G)
$$
and $V_2$ is an independent set.
\end{Thm}
\begin{proof}
Suppose that $l=l(G)$ and
\begin{align*}
&L_1=x_{1,1}-x_{1,2}-\cdots-x_{1,l}\\
&L_2=x_{2,1}-x_{2,2}-\cdots-x_{2,l}\\
&\qquad\cdots\cdots\\
&L_t=x_{t,1}-x_{t,2}-\cdots-x_{t,l}
\end{align*}
are all paths of $G$ with the length $l$. Below we shall construct
some sets $U_i$ and $I_i$. Let $U_1=\{x_{1,1}\}$ and
$$
I_1=\{1\leq i\leq t:\, U_1\cap L_i=\emptyset\}.
$$
For $j\geq 2$,  if $I_{j-1}\not=\emptyset$, then let $i'=\min
I_{j-1}$, $U_j=U_{j-1}\cup\{x_{i',1}\}$ and
$$
I_j=\{1\leq i\leq t:\, U_j\cap L_i=\emptyset\}.
$$
Of course, if $I_{j-1}=\emptyset$, then stop this process. Suppose
that we finally get the vertex set $U_s$. Assume that
$U_s=\{x_{i_1,1},x_{i_2,1},\ldots,x_{i_s,1}\}$ where
$1=i_1<i_2<\cdots<i_s$. Let $V_2=U_s$ and $V_1=V\setminus V_2$.
First, we claim that $V_2$ is an independent set. Assume on the
contrary that there exist $1\leq a<b\leq s$ such that $x_{i_a,1}$
and $x_{i_b,1}$ are adjacent in $G$. By the construction of $U_s$,
we have $x_{i_a,1}$ doesn't lie in the path $L_{i_b}$. Thus
$$
x_{i_a,1}-x_{i_b,1}-x_{i_b,2}-\cdots-x_{i_b,l}
$$
forms a path with the length $l+1$. It is impossible since $l(G)=l$.
Second, by noting that $I_s=\emptyset$, we have $V_2\cap
L_i\not=\emptyset$ for any $1\leq i\leq t$. Hence $l(G[V_1])<l$
since $L_1,\ldots,L_t$ are all paths of $G$ with the length $l$.
\end{proof}

\begin{Ack}
The authors thank Professor Zhi-Wei Sun for his useful suggestions.
And the first author thanks Professor Yu-Sheng Li for his helpful
discussions.
\end{Ack}

\end{document}